# An elementary proof that P ≠ NP[1]

Bhupinder Singh Anand[2]

We show that, if PA has no non-standard models, then P ≠ NP. We then give an elementary proof that PA has no non-standard models.

## 1. Introduction

In a paper presented to ICM 2002, Ran Raz comments [Ra02]:

> "A Boolean formula $f(x_1, ..., x_n)$ is a tautology if $f(x_1, ..., x_n) = 1$ for every $x_1, ..., x_n$. A Boolean formula $f(x_1, ..., x_n)$ is unsatisfiable if $f(x_1, ..., x_n) = 0$ for every $x_1, ..., x_n$. Obviously, $f$ is a tautology if and only if $\sim f$ is unsatisfiable.
>
> Given a formula $f(x_1, ..., x_n)$, one can decide whether or not $f$ is a tautology by checking all the possibilities for assignments to $x_1, ..., x_n$. However, the time needed for this procedure is exponential in the number of variables, and hence may be exponential in the length of the formula $f$. …
>
> P ≠ NP is the central open problem in complexity theory and one of the most important open problems in mathematics today. The problem has thousands of equivalent formulations. One of these formulations is the following:

---

Is there a polynomial time algorithm A that gets as input a Boolean formula *f* and outputs 1 if and only if *f* is a tautology? P ≠ NP states that there is no such algorithm."

## 2. Gödel's proof and the PvNP problem

We note, first, that, in his seminal paper on undecidable arithmetical propositions [Go31], Gödel has defined a formula, $[R(x)]$[3], and shown that:

(*i*)    $[R(x)]$ is constructible in a standard, first-order, Peano Arithmetic, PA;

(*ii*)   we can prove, meta-mathematically, that $[R(x)]$ translates as an arithmetical tautology, $R(x)$, under the standard interpretation (*Appendix C*), M, of the Arithmetic;

(*iii*)  $[R(x)]$ is not provable in the Arithmetic.

We consider, next:

**Definition 1**: A total number-theoretical relation, $R(x_1, x_2, ..., x_n)$, when treated as a Boolean function, is Turing-decidable in M if, and only if, it is instantiationally equivalent to a number-theoretic relation, $S(x_1, x_2, ..., x_n)$, and there is a Turing-machine T such that, for any given natural number sequence, $(a_1, a_2, ..., a_n)$, T will compute $S(a_1, a_2, ..., a_n)$ as either TRUE, or as FALSE.

**Definition 2**: A total number-theoretical relation, $R(x_1, x_2, ..., x_n)$, when treated as a Boolean function, is Turing-computable in M if, and only if, there is a Turing-machine T such that, for any given natural number sequence, $(a_1, a_2, ..., a_n)$, T will compute $R(a_1, a_2, ..., a_n)$ as either TRUE, or as FALSE.

---

[3] We use square brackets to indicate that the expression inside the brackets is to be viewed purely as a formula (*i.e., as an uninterpreted, syntactical, string of symbols*) of a formal system. Absence of brackets indicates that the expression is to be viewed as a specific interpretation of the formula (*under Tarski's definitions of the satisfiability and truth of a formula of a formal system under a given interpretation*).

Now, in [Go31], Gödel's reasoning only shows that $R(x)$ is Turing-decidable as always TRUE (*when treated as a Boolean function*).

However, the question remains:

Is $R(x)$ also Turing-computable as always TRUE (*when treated as a Boolean function*)?

The distinction[4] assumes significance in the light of the PvNP problem [Cook].

For, if we assume, first, that every total arithmetical relation, which is Turing-computable as always TRUE, is the standard interpretation of a PA-provable formula, then $R(x)$ is not Turing-computable as always TRUE, and, so, $P \neq NP$.

If we assume, however, that there is a total arithmetical relation that is Turing-computable as always TRUE, but which is not the standard interpretation[5] of a PA-provable formula, then this implies that there is a non-standard model[6] of PA.

---

[4] Note that, although decidability of a number-theoretic relation, which can be interpreted as a Turing-algorithm in the classical sense (cf. [Me64], p229-237), necessarily implies Turing-computability, the converse need not hold.

Thus, classically, such a number-theoretic relation, $R(x)$, would be decidable as true if it is instantiationally equivalent, for any given natural number, to a number-theoretic relation, $S(x)$, which can be interpreted as a Turing-algorithm in the classical sense, and which is Turing-computable as TRUE for all natural numbers.

However, we cannot conclude from this, without proof, that $R(x)$ must also be interpretable as a Turing-algorithm in the classical sense, and be Turing-computable as TRUE for all natural numbers.

In fact the PvsNP problem, as enunciated by Raz [Ra02], is equivalent to asking whether there are number-theoretical relations that are Turing-decidable as always TRUE, but not Turing-computable as always TRUE.

[5] The word "interpretation" may be used both in its familiar, linguistic, sense, and in a mathematically precise sense; the appropriate meaning is usually obvious from the context.

We conclude that, if PA has no non-standard models, then, under the above expression [Ra02] of the PvNP problem, $P \neq NP$.

## 3. Standard, first-order, PA has no non-standard model

We now give an elementary proof that, if PA is a standard, first-order, Peano Arithmetic - as defined in Appendix A and Appendix B - then PA has no non-standard model.

We denote by $G(x)$ the PA-formula:

   $[x=0 \text{ v } \sim(\forall y)\sim(x=y')]$.

This translates, under every interpretation of PA, as:

   Either $x$ is 0, or $x$ is a 'successor'.

Now, in every interpretation of PA, we have that:

   (*a*) $G(0)$ is true;

---

Mathematically, following Tarski ([Me64], §2, p49): "An *interpretation* consists of a non-empty set D, called the *domain* of the interpretation, and an assignment to each predicate letter $A_j^n$ of an *n*-place relation in D, to each function letter $f_j^n$ of an *n*-place operation in D (*i.e., a function from $D^n$ into D*), and to each individual constant $a_i$ of some fixed element of D. Given such an interpretation, variables are thought of as ranging over the set D, and $\sim$, $\Rightarrow$, and quantifiers are given their usual meaning. (*Remember that an n-place relation in D can be thought of as a subset of $D^n$, the set of all n-tuples of elements of D.*)"

We note that the interpreted relation $R'(x)$ is obtained from the formula $[R(x)]$ of a formal system P by replacing every primitive, undefined symbol of P in the formula $[R(x)]$ by an interpreted mathematical symbol. So the P-formula $[(\forall x)R(x)]$ interprets as the sentence $(\forall x)R'(x)$, and the P-formula $[\sim(\forall x)R(x)]$ as the sentence $\sim(\forall x)R'(x)$.

[6] We define [cf. [Me64], p49-53] a model as an interpretation of a set of well-formed formulas if, and only if, every well-formed formula of the set is true for the interpretation (*by Tarski's definitions of the satisfiability and truth of the formulas of a formal language under an interpretation*).

We define a non-standard model, say M', of first-order Peano Arithmetic as one in which, if $[(\forall x)R(x)]$ is Gödel's undecidable formula, then $\sim R(s)$ holds for some $s$ in the domain of M' that is not a natural number (*i.e., s is not a successor of 0*).

(*b*) If *G*(*x*) is true, then *G*(*x*') is true.

It follows, from Gödel's completeness theorem[7], that:

(*a*) [*G*(0)] is provable in PA;

(*b*) [*G*(*x*) => *G*(*x*')] is provable in PA.

We also have, by Generalisation (*Appendix A*), that:

(*c*) [($\forall x$)(*G*(*x*) => *G*(*x*'))] is provable in PA;

From the Induction axiom S9 (*Appendix B*), we thus have that:

(*d*) [($\forall x$)*G*(*x*)] is provable in PA.

We conclude that, except 0, every element in the domain of any interpretation of PA is a 'successor' of 0.

## 5. Conclusions

Since, by definition, the 'successors' of 0 are the natural numbers, it follows that:

**Theorem 1**: There are no non-standard models of a standard, first-order, Peano Arithmetic.

**Corollary 1**: $P \neq NP$.

**Corollary 2**: An arithmetical relation is Turing-computable as always TRUE if, and only if, it is the standard interpretation of a PA-provable formula.

---

[7] Cf. [Me64], p68.

## Appendix A: Logical Symbols, Axioms, and Rules of Inference of PA

The primitive logical symbols of PA are:

| $\rightarrow$ | & | v | ~ | $\forall$ | $x, y, ...$ | $a, b, ...$ |
|---|---|---|---|---|---|---|
| implies | and | or | not | for all | variables | constant terms |

If $A$, $B$, $C$ are well-formed formulas of PA, then the logical axioms of PA are:

(1) $A \rightarrow (B \rightarrow A)$;

(2) $(A \rightarrow (B \rightarrow C)) \rightarrow ((A \rightarrow B) \rightarrow (A \rightarrow C))$;

(3) $(\sim B \rightarrow \sim A) \rightarrow ((\sim B \rightarrow A) \rightarrow B)$;

(4) $(\forall x)A(x) \rightarrow A(t)$, if $A(x)$ is a well-formed formula of PA, and $t$ is a term of PA free for $x$ in $A(x)$;

(5) $(\forall x)(A \rightarrow B) \rightarrow (A \rightarrow (\forall x)B)$, if $A$ is a well-formed formula of PA containing no free occurrences of $x$.

The rules of inference of PA are:

(*i*)   Modus Ponens: $B$ follows from $A$ and $A \rightarrow B$;

(*ii*)  Generalisation: $(\forall x)A$ follows from $A$.

By Tarski's definitions of the satisfiability and truth of the formulas of PA under an interpretation, when the rules of inference are applied to true well-formed formulas of PA under a given interpretation, then the results of these applications are also true (*i.e., every theorem of PA is true in any model of PA*). (cf. [Me64], p57)

## Appendix B: Primitive Symbols and Proper Axioms of PA

In a standard, first order, Peano Arithmetic, ([Me64], p103):

(*a*)  PA has a single predicate letter, $A^2_1$ (as usual, we write "$t = s$" for $A^2_1(t, s)$;

(*b*)  PA has one individual constant $a_1$ (written, as usual, " 0 ");

(*c*)  PA has three function letters $f^1_1, f^2_1, f^2_2$. We shall write "$t'$" instead of $f^1_1(t)$; "$t+s$" instead of $f^2_1(t, s)$; and "$t*s$" instead of $f^2_2(t, s)$.

The proper axioms of PA are:

(S1)  $(x_1 = x_2) \rightarrow ((x_1 = x_3) \rightarrow (x_2 = x_3))$;

(S2)  $(x_1 = x_2) \rightarrow (x_1' = x_2')$;

(S3)  $0 \neq (x_1)'$;

(S4)  $((x_1)' = (x_2)') \rightarrow (x_1 = x_2)$;

(S5)  $(x_1 + 0) = x_1$;

(S6)  $(x_1 + x_2') = (x_1 + x_2)'$;

(S7)  $(x_1 * 0) = 0$

(S8)  $(x_1 * (x_2')) = ((x_1 * x_2) + x_1)$;

(S9)  For any well-formed formula $F(x)$ of PA:
$F(0) \rightarrow ((\forall x)(F(x) \rightarrow F(x')) \rightarrow (\forall x)F(x))$.

## Appendix C: The standard interpretation, M, of PA ([Me64], p107)

The standard interpretation, M, of PA, is taken to be the intuitive arithmetic of the natural numbers, as expressed by Dedekind's semi-axiomatic formulation of the Peano Postulates, where:

(*i*) the integer 0 is the interpretation of the PA-symbol " 0 ";

(*ii*) the successor operation (addition of 1) is the interpretation of the PA-symbol " ′ ";

(*iii*) ordinary addition and multiplication are the interpretations of the PA-symbols " + " and " * ";

(*iv*) the interpretation of the PA-symbol " = " is the identity relation.